\documentclass[11pt]{article}
\usepackage{amsmath,amsthm,amsfonts,amssymb,amscd, amsxtra}

\newcommand{\banacha}{X}
\newcommand{\banachb}{Y}
\newtheorem{theorem}{Theorem}
\newtheorem{lemma}[theorem]{Lemma}

\newtheorem{corollary}[theorem]{Corollary}
\newtheorem{proposition}[theorem]{Proposition}

\newtheorem{remark}{Remark}
\newtheorem{example}{Example}
\begin{document}
\title{Local convergence of Newton's  method  under   majorant condition}

\author{ O. P. Ferreira\thanks{IME/UFG,  CP-131, CEP 74001-970 - Goi\^ania, GO, Brazil ({\tt
      orizon@mat.ufg.br}). The author was supported  by
     CNPq Grants 302618/2005-8 and 475647/2006-8,  PRONEX--Optimization(FAPERJ/CNPq) and FUNAPE/UFG.}  
}
\date{May 13, 2009}
\maketitle
\begin{abstract}
A local convergence analysis of Newton's method for solving nonlinear equations, under a majorant condition, is presented in this paper. Without assuming convexity of the derivative of the majorant function, which relaxes the Lipschitz condition on the operator under consideration,  convergence, the biggest range for uniqueness of the solution,  the optimal convergence radius and results on the convergence rate are established. Besides, two special cases of the general theory are presented as an application.


\end{abstract}
\section{Introduction}\label{sec:int}

The Newton's  method and its variant are  powerful tools for solving nonlinear equation in real or complex Banach space. In the last few years, a couple of papers have dealt with the issue of local and semi-local convergence analysis of Newton's method and its variants  by relaxing the assumption of Lipschitz continuity of the derivative of the function, which define  the nonlinear equation in consideration, see \cite{C08}, \cite{C10}, \cite{Cl08}, \cite{F08}, \cite{FG09},  \cite{FS09},  \cite{Hu04}, \cite{W00}, \cite{W03}.

In  \cite{F08} and  \cite{W00}, under a majorant condition and generalized Lipschitz condition, respectively,   local convergence, quadratic rate and estimate of the best possible convergence radius of the Newton's method were established, as well as uniqueness of solution for the nonlinear equation in question. The analysis presented in   \cite{F08}, \emph{convexity} of the derivative of the scalar majorant function was assumed and  in  \cite{W00} the {\it nondecrement} of the positive integrable function which defines the generalized Lipschitz condition.  These assumptions seem  to be actually natural in local analysis of Newton's method. Even though,  convergence, uniqueness, superlinear  rate and estimate of the best possible convergence radius will be established in this paper without assuming convexity of the derivative of the majorant function or that the function which defines the generalized Lipschitz condition  is  nondecreasing. In particular, this analysis shows  that the  convexity of the derivative of the majorant function or that the function which defines the generalized Lipschitz condition  is  nondecreasing are  needed only to obtain quadratic convergence rate of the sequence generated by Newton's Method. Also, as in \cite{F08}, the analysis presented provides a clear relationship between the majorant function  with the nonlinear operator under consideration. Besides improving the convergence theory  this analysis permits  us to obtain two new important special cases, namely, \cite{Hu04} and \cite{W03} (see also, \cite{F08} and \cite{W00}) as an application. It is worth pointing out that the majorant condition used here is equivalent to Wang's condition (see  \cite{W03}) always the derivative of majorant function is convex.

The organization of the paper is as follows.
In Section \ref{sec:int.1},  some notations and one basic result used in the paper are presented. In Section \ref{lkant}, the main result is stated and  in   Section \ref{sec:PR} some properties of the majorant function are established and the main relationships  between the majorant function and the nonlinear operator used in the paper are presented. In Section \ref{sec:UnOpBall}, the uniqueness of the solution and the optimal convergence radius are obtained. In  Section \ref{sec:proof} the main result is proved and two applications of this result are given in Section \ref{apl}. Some final remarks are made in Section~\ref{fr}.
\subsection{Notation and auxiliary results} \label{sec:int.1}
The following notations and results are used throughout our presentation. Let $\banacha$,\, $\banachb$ be  Banach spaces. The open and closed ball at $x$ are denoted, respectively, by
$$
B(x,\delta) = \{ y\in X ;\; \|x-y\|<\delta \}\;\;\; \mbox{and}\;\;\;B[x,\delta] = \{ y\in X ;\; \|x-y\|\leqslant \delta
\}.
$$
Let $\Omega\subseteq \banacha$ an open set. The Fr\'echet derivative of $F:{\Omega}\to \banachb$ is the linear map $F'(x):\banacha \to \banachb$.
\begin{lemma}[Banach's Lemma] \label{lem:ban}
Let $B:\banacha \to \banacha$ be  bounded linear operator. If   $I:\banacha \to \banacha$ is the identity operator and $\|B-I\|<1$,  then $B$ is invertible and 
$
\|B^{-1}\|\leq 1/\left(1-  \|B-I\|\right).
$
\end{lemma}
\section{Local analysis  for  Newton's method } \label{lkant}
Our goal is to state and prove a local theorem for Newton's method,  which generalize Theorem~2.1 of \cite{F08}. First, we will
 prove some results regarding the scalar majorant function, which relaxes the Lipschitz condition. Then we will establish the main relationships between the majorant function and the nonlinear function. We will also prove the uniqueness of the solution in a suitable region and the optimal ball of convergence. Finally, we will show well definedness of Newton's method  and convergence, also results on the convergence rates will be given.  The statement of the theorem~is:
\begin{theorem}\label{th:nt}
Let $\banacha$, $\banachb$ be Banach spaces, $\Omega\subseteq \banacha$ an open set and
  $F:{\Omega}\to \banachb$  a continuously
  differentiable function. 
 Let $x_*\in \Omega$, $R>0$  and  $\kappa:=\sup\{t\in [0, R): B(x_*, t)\subset \Omega\}$.
  Suppose that $F(x_*)=0$,  $F '(x_*)$ is  invertible and there exist an $f:[0,\; R)\to \mathbb{R}$ continuously differentiable such~that  
  \begin{equation}\label{Hyp:MH}
\left\|F'(x_*)^{-1}\left[F'(x)-F'(x_*+\tau(x-x_*))\right]\right\| \leq    f'\left(\|x-x_*\|\right)-f'\left(\tau\|x-x_*\|\right),
  \end{equation}
  for all $\tau \in [0,1]$, $x\in B(x_*, \kappa)$ and 
\begin{itemize}
  \item[{\bf h1)}]  $f(0)=0$ and $f'(0)=-1$;
  \item[{\bf  h2)}]  $f'$ is  strictly increasing.
\end{itemize}
  Let  $\nu:=\sup\{t\in [0, R): f'(t)< 0\},$ $\rho:=\sup\{\delta\in (0, \nu): [f(t)/f'(t)-t]/t<1,\, t\in (0, \delta)\} $  and 
  $$
r:=\min \left\{\kappa, \,\rho \right\}.
  $$
Then the sequences  with starting points $x_0\in B(x_*, r)/\{x_*\}$ and $t_0=\|x_0-x_*\|$, respectively, namely  
\begin{equation} \label{eq:DNS}
    x_{k+1} ={x_k}-F'(x_k) ^{-1}F(x_k), \qquad   t_{k+1} =|t_k-f(t_k)/f'(t_k)|, \qquad k=0,1,\ldots\,, 
\end{equation}
are well defined; $\{t_k\}$ is strictly decreasing, is contained in $(0, r)$ and converges to $0$ and $\{x_k\}$ is contained in $B(x_*,r)$ and  converges to the point $x_*$ which is the unique zero of $F$ in $B(x_*, \sigma)$, where $\sigma:=\sup\{t\in (0, \kappa): f(t)< 0\}$ and there hold:
\begin{equation}
    \label{eq:q2}
    \lim_{k \to \infty}\left[\|x_{k+1}-x_*\|\big{/}\|x_k-x_*\|\right]=0, \qquad \lim_{k \to \infty}[t_{k+1}/t_k]=0. 
  \end{equation}
Moreover, if  $f(\rho)/(\rho f'(\rho)-1=1$ and $\rho<\kappa$ then $r=\rho$ is the best  possible convergence radius.

\noindent
If, additionally, given $0\leq p\leq1$
\begin{itemize}
  \item[{\bf  h3)}] the function  $(0,\, \nu) \ni t \mapsto [f(t)/f'(t)-t]/t^{p+1}$  is  strictly increasing,
\end{itemize}
then the sequence $\{t_{k+1}/t_k^{p+1}\}$ is strictly decreasing  and there holds
\begin{equation}
    \label{eq:q3}
\|x_{k+1}-x_*\| \leq \big[t_{k+1}/t_k^{p+1}\big]\,\|x_k-x_*\|^{p+1}, \qquad k=0,1,\ldots\,. 
  \end{equation}
\end{theorem}
\begin{remark} \label{r:rqc}
The first equation  in \eqref{eq:q2} means that  $\{x_k\}$ converges  superlinearly  to $x_*$. Moreover, because  the sequence $\{t_{k+1}/t_k^{p+1}\}$ is strictly decreasing  then
$
t_{k+1}/ t_k^{p+1}\leq t_1/t_0^{p+1},
$
for $k=0,1,  \ldots$. 
So, the inequality in \eqref{eq:q3} implies $\|x_{k+1}-x_*\| \leq [t_{1}/t_0^{p+1}]\|x_k-x_*\|^{p+1},$
for $ k=0,1,\ldots\, $. As a consequence,  if $p=0$ then $\|x_{k}-x_*\| \leq t_0[t_1/t_0]^k$ for $k=0,1,\ldots\,$ and if $0<p\leq 1$ then
$$
\|x_{k}-x_*\| \leq t_0\left(t_1/t_0\right)^{[(p+1)^k-1]/p}, \qquad k=0,1,\ldots\,.
$$
\end{remark}
\begin{example} \label{ex:mf}
The following  continuously differentiable functions satisfy {\bf h1}, {\bf h2} and  {\bf h3}:
\begin{itemize}
  \item[{i)}] $f: [0, +\infty)\to \mathbb{R}$ such that $f(t)=t^{1+p}-t $;
  \item[{ii)}] $f: [0, +\infty)\to \mathbb{R}$ such that $f(t)=\mbox{e}^{-t}+t^2-1$.
\end{itemize}
Letting $0<p<1$,  the derivative of first function is not convex, as well as of the second.
\end{example}
Similarly to the proof of  Proposition $2.6$ in  \cite{F08}, always  $f$ has derivative $f'$ convex,  we can prove that  {\bf  h3} holds with $p=1$.  In this case,   the Newton's sequence  converges  with quadratic rate. Indeed,  the next example  shows that the convexity of $f'$ was necessary in  \cite{F08} to obtain quadratic  convergence rate.
\begin{example}
Let $g:  \mathbb{R}\to \mathbb{R}$ be given by $g(t)=t^{5/3}-t.$ Note that $g(0)=0$,  $g'(0)=-1$ and letting $p=2/3$ in   Example \ref{ex:mf} the function $f$ is  a mojorant function to $g$. The Newton's method applied to $g$ with starting point $t_0$ ``near'' $0$ generates the following sequence:
$$
t_{k+1}=\left(2\,t_{k}^{5/3}\right)\big{/}\left(5\,t_{k}^{2/3}-3\right), \qquad k=0,1,\ldots\,.
$$
Theorem \ref{th:nt} implies that the sequence  $\{t_k\}$ converges to $0$ with superlinear rate. It easy to see that $\{t_k\}$ does not converges to $0$ with quadratic rate. So,  in particular, it follows from \cite{F08} that there is none majorant function having convex derivative for the function $g$.
\end{example}
From now on, we assume that the hypotheses of Theorem \ref{th:nt}
hold, with the exception of {\bf h3} which will be considered to
hold only when explicitly stated.

\subsection{Preliminary results} \label{sec:PR}
In this section, we will prove all  statements in Theorem~\ref{th:nt} regarding  the sequence $\{t_k\}$ associated to the  majorant function. The main relationships between the  majorant function  and the nonlinear operator will be also established, as well as the results in Theorem~\ref{th:nt}  related to the uniqueness of the solution and the optimal convergence radius.
\subsubsection{The scalar sequence} \label{sec:PMF}
In this section,  we will prove the statements in Theorem~\ref{th:nt} involving  $\{t_k\}$.
First, we will prove that the constants $\kappa$,  $\nu$,  $\rho$ and $\sigma$ are positive. We beginning proving that    $\kappa$, $\nu$ and $\sigma$ are positive.
\begin{proposition}  \label{pr:incr1}
The constants $ \kappa,\, \nu $ and $\sigma$ are positive and $t-f(t)/f'(t)<0,$ for all $t\in (0,\,\nu).$
\end{proposition}
\begin{proof}
Since $\Omega$ is open and $x_*\in \Omega$, we can immediately conclude that $\kappa>0$. As $f'$ is continuous in $0$ with  $f'(0)=-1$, there exists $\delta>0$ such that $f'(t)<0$ for all $t\in (0,\, \delta).$ So,  $\nu>0$. Now, because $f(0)=0$ and  $f'(0)=-1$,  there exists $\delta>0$ such that $f(t)<0$ for all $t\in (0, \delta)$. Hence $\sigma>0$. 

It remains to show that $t-f(t)/f'(t)<0,$ for all $t\in (0,\,\nu).$ Since $f'$ is strictly increasing,  $f$ is strictly convex. So, 
$
0=f(0)>f(t)-tf'(t),
$
for $t\in  (0,\, R).$
If $t\in (0, \,\nu)$ then $f'(t)<0$, which, combined with last inequality yields the desired inequality.
\end{proof}

According to {\bf h2} and  definition of $\nu$, we have  $f'(t)< 0$ for all
$t\in[0, \,\nu)$.  Therefore, the Newton iteration map for  $f$ is well defined in
$[0,\, \nu)$. Let us call it $n_f$:
\begin{equation} \label{eq:def.nf}
  \begin{array}{rcl}
  n_f:[0,\, \nu)&\to& (-\infty, \, 0]\\
    t&\mapsto& t-f(t)/f'(t).
  \end{array}
\end{equation}
\begin{proposition}  \label{pr:incr3}
$
\lim_{t\to 0}|n_f(t)|/t=0.
$
As a consequence,  $\rho>0 $ and 
$|n_f(t)|<t$ for all $ t\in (0, \, \rho)$.
\end{proposition}
\begin{proof}
Using definition \eqref{eq:def.nf},  Proposition \ref{pr:incr1},  $f(0)=0$, and definition of $\nu$,  simple algebraic manipulation gives 
\begin{equation} \label{eq:rho}
\frac{|n_f(t)|}{t}= [f(t)/f'(t)-t]/t=\frac{1}{f'(t)} \frac{f(t)-f(0)}{t-0}-1, \qquad t\in (0,\,\nu).
\end{equation}
Because  $f'(0)=-1\neq 0$ the first statement follows by taking limit in~\eqref{eq:rho}, as $t$ goes to $0$.

Since $\lim_{t\to 0}|n_f(t)|/t=0$,  first equality in  \eqref{eq:rho} implies that there exists $\delta>0$ such that 
$$
0<[f(t)/f'(t)-t]/t<1, \qquad  \;  t\in (0, \delta). 
$$
So,  we conclude  that $\rho$  is positive. Therefore,  the first equality in \eqref{eq:rho} together  definition of  $\rho$  implies that $|n_f(t)|/t=[f(t)/f'(t)-t]/t<1,$ for  all $t\in (0, \rho)$, as required.
\end{proof}

Using \eqref{eq:def.nf}, it easy to see that  the sequence $\{t_k \}$ is equivalently defined as
\begin{equation} \label{eq:tknk}
 t_0=\|x_0-x_*\|, \qquad t_{k+1}=|n_f(t_k)|, \qquad k=0,1,\ldots\, .
\end{equation}
\begin{corollary} \label{cr:kanttk} 
The sequence $\{t_k\}$ is well defined, is strictly decreasing and is contained in $(0, \rho)$. Moreover,  $\{t_k\}$ converges to $0$ with superlinear rate, i.e., 
$
\lim_{k\to \infty}t_{k+1}/t_k=0.
$
If, additionally, {\bf  h3} holds then the sequence $\{t_{k+1}/t_k^{p+1}\}$ is strictly decreasing.
\end{corollary}
\begin{proof}
Since $0<t_0=\|x_0-x_*\|<r\leq\rho$, using  Proposition~\ref{pr:incr3} and \eqref{eq:tknk}  it is simple to conclude that $\{t_k\}$ is well defined, is strictly decreasing and is contained in $(0, \rho)$. So, we have proved the first statement of the corollary.

Because  $\{t_k \}\subset (0, \rho)$ is  strictly decreasing  it converges. So,  $\lim_{k\to \infty}t_{k}=t_*$ with $0\leq t_*<\rho$ which together with \eqref{eq:tknk} implies $0\leq t_{*}=|n_f(t_*)|$. But, if $t_*\neq 0$ then Proposition~\eqref{pr:incr3}  implies $|n_f(t_*)|<t_*$, hence $t_*=0$. Now, as $\lim_{k\to \infty}t_{k}=0$. Thus, definition of $\{t_k\}$ in  \eqref{eq:tknk} and first statement in Proposition~\ref{pr:incr3}  imply that  $\lim_{k\to \infty}t_{k+1}/t_k=\lim_{k\to \infty}|n_f(t_k)|/t_k=0$  and the  second statement is proved.

   Since $\{t_k \}$ is  strictly decreasing, the last statement is an immediate consequence of {\bf  h3}.
\end{proof}

\subsubsection{Relationship between the majorant function and the nonlinear operator} \label{sec:MFNLO}
In this section, we will present the main  relationships between the majorant function $f$ and the nonlinear operator $F$. 
\begin{lemma} \label{wdns}
If \,\,$\| x-x_*\|<\min\{\kappa, \nu\}$,  then $F'(x) $ is invertible and
$$
\|F'(x)^{-1}F'(x_*)\|\leqslant  1/|f'(\| x-x_*\|)|. 
$$
In particular, $F'$ is invertible in $B(x_*, r)$.
\end{lemma}
\begin{proof}
For proving this lemma is not necessary the assumption {\bf  h3} neither that the derivative of the majorant function is convex. The proof follows the same pattern of Lemma~2.9 of \cite{F08}.
\end{proof}

Newton iteration at a point happens to be a zero of the linearization
of $F$ at such a point.  So, we study the  linearization error  at a point
in $\Omega$
\begin{equation}\label{eq:def.er}
  E_F(x,y):= F(y)-\left[ F(x)+F'(x)(y-x)\right],\qquad y,\, x\in \Omega.
\end{equation}
We will bound this error by the error in the linearization on the
majorant function $f$
\begin{equation}\label{eq:def.erf}
        e_f(t,u):= f(u)-\left[ f(t)+f'(t)(u-t)\right],\qquad t,\,u \in [0,R).
\end{equation}
\begin{lemma}  \label{pr:taylor}
If  $\|x_*-x\|< \kappa$, then  
$
\|F'(x_*)^{-1}E_F(x, x_*)\|\leq e_f(\|x-x_*\|, 0). 
$
\end{lemma}
\begin{proof}
For proving this lemma is not necessary the assumption {\bf  h3} neither that the derivative of the majorant function is convex. The proof follows the same pattern of Lemma~2.10 of \cite{F08}.
\end{proof}

Lemma \ref{wdns} guarantees, in particular,  that  $F'$ is invertible in $B(x_*, r)$ and consequently, the Newton iteration map is well-defined.  Let us call $N_{F}$, the Newton
iteration map for $F$ in that region:
\begin{equation} \label{NF}
  \begin{array}{rcl}
  N_{F}:B(x_*, r) &\to& \banachb\\
    x&\mapsto& x-F'(x)^{-1}F(x).
  \end{array}
\end{equation}
Now,  we establish an important  relationship between  the Newton iteration maps $n_{f}$ and $ N_{F}$. As a consequence,   we obtain that $B(x_*, r)$ is invariant under  $ N_{F}$. This result will be important to assert the well definition of the Newton method.
\begin{lemma} \label{le:cl}
If   $\|x-x_*\|< r$  then 
$
\|N_F(x)-x_*\|\leq  |n_f(\|x-x_*\|)|.
$
As a consequence,  $$N_{F}(B(x_*, r))\subset B(x_*, r).$$
\end{lemma}
\begin{proof}
Since $F(x_*)=0$, the inequality is trivial for $x=x_*$. Now assume that  $0<\|x-x_*\|\leq t$. 
Lemma \ref{wdns} implies that  $F'(x) $ is invertible. Thus, because $F(x_*)=0$, direct  manipulation yields  
$$
x_*-N_F(x)=-F'(x)^{-1}\left[ F(x_*)-F(x)-F'(x)(x_*-x)\right]= -F'(x)^{-1}E_F(x,x_*).
$$
Using  the above equation,  Lemma \ref{wdns} and Lemma \ref{pr:taylor}, it is easy to conclude that
$$
 \|x_*-N_F(x)\|\leq\| -F'(x)^{-1}F'(x_*)\|\|  F'(x_*)^{-1}E_F(x,x_*)\|\leq e_f(\|x-x_*\|, 0)/|f'(\|x-x_*\|)|.
$$
On the other hand, taking into account  that $f(0)=0$,   the definitions of $e_f$ and  $n_f$  imply that
$$
e_f(\|x-x_*\|, 0)/|f'(\|x-x_*\|)|=|n_f(\|x-x_*\|)|.
$$
So,  the first statement follows by combining two above expressions.

Take $x\in B(x_*, r)$.  Since $\|x-x_*\|<r$ and $ r\leq \rho$, the first part  together with the second part of Proposition~\ref{pr:incr3} imply that 
$
\|N_F(x)-x_*\|\leq |n_f(\|x-x_*\|)|<\|x-x_*\|
$
and the  last result  follows.
\end{proof}
\begin{lemma} \label{le:cl2}
If {\bf h3} holds and $\|x-x_*\|\leq t<r$ then 
$
\|N_F(x)-x_*\|\leq [ |n_f(t)|/t^{p+1}]\,\|x-x_*\|^{p+1}.
$
\end{lemma}
\begin{proof}
The inequality is trivial for $x=x_*$. 
If $0<\|x-x_*\|\leq t$ then  assumption {\bf  h3} and   \eqref{eq:def.nf}  give $|n_f(\|x-x_*\|)|/\|x-x_*\|^{p+1}\leq |n_f(t)|/t^{p+1}$. So, using   Lemma~\ref{le:cl} the statement follows.
\end{proof}
\subsection{Uniqueness and  optimal convergence radius} \label{sec:UnOpBall}
In this section we will obtain the uniqueness of the solution and the optimal convergence radius.
\begin{lemma} \label{pr:uniq}
 The point  $x_*$ is the unique zero of $F$ in $B(x_*, \sigma)$.
\end{lemma}
\begin{proof}
For proving this lemma is not necessary the assumption {\bf  h3} neither that the derivative of the majorant function is convex. The proof follows the same pattern of Lemma~2.13 of \cite{F08}.
\end{proof}
\begin{lemma} \label{pr:best}
If  $f(\rho)/(\rho f'(\rho))-1=1$ and $\rho < \kappa$, then  $r=\rho$ is the optimal convergence radius.
\end{lemma}
\begin{proof}
The proof follows the same pattern of Lemma~2.15 of \cite{F08}.
\end{proof}
\subsection{The Newton's sequence} \label{sec:proof}
In this section,  we will prove the statements in Theorem~\ref{th:nt} involving  the Newton's sequence $\{x_k\}$. First, note that the first equation in \eqref{eq:DNS} together with \eqref{NF} implies that   the sequence $\{x_k\}$  satisfies 
\begin{equation} \label{NFS}
x_{k+1}=N_F(x_k),\qquad k=0,1,\ldots \,,
\end{equation}
which is indeed an equivalent definition of this sequence.
\begin{proposition}\label{pr:nthe}
The sequence $\{x_k\}$ is well defined, is contained in $B(x_*,r)$ and  converges to the point $x_*$ the unique zero of $F$ in $B(x_*, \sigma)$ and there hold:
\begin{equation} \label{eq:q2e}
    \lim_{k \to \infty}\left[\|x_{k+1}-x_*\|\big{/}\|x_k-x_*\|\right]=0. 
  \end{equation}
If, additionally, {\bf  h3} holds then the sequences  $\{x_k\}$  and $\{t_k\}$ satisfy
\begin{equation} \label{eq:q3e}
\|x_{k+1}-x_*\| \leq \big[t_{k+1}/t_k^{p+1}\big]\,\|x_k-x_*\|^{p+1}, \qquad k=0,1,\ldots\,. 
  \end{equation}
\end{proposition}
\begin{proof} 
As $x_0\in B(x_*,r)$ and $r\leq \nu$, combining  \eqref{NFS},  inclusion  $N_{F}(B(x_*, r)) \subset B(x_*, r)$ in Lemma~\ref{le:cl} and Lemma~\ref{wdns}, it is easy to conclude that   $\{x_k\}$  is well defined and remains  in $B(x_*,r)$. 

We are going to prove that $\{x_k \}$ converges towards $x_*$. Since $\|x_{k}-x_*\|<r\leq \rho$,  for $ k=0,1,\ldots \,$,  we obtain from  \eqref{NFS}, Proposition~\ref{le:cl} and Proposition~\ref{pr:incr3} that
\begin{equation}\label{eq:conv1}
\|x_{k+1}-x_*\|=\|N_F(x_k)-x_*\|\leq |n_f(\|x_{k}-x_*\|)|<\|x_{k}-x_*\|,\qquad  k=0,1,\ldots \,.
\end{equation}
So, $\{\|x_{k}-x_*\| \}$ is strictly decreasing and convergent. Let $\ell_*=\lim_{k\to \infty}\|x_{k}-x_*\|$. Because  $\{\|x_{k}-x_*\| \}$ rest in $(0, \,\rho)$ and is strictly decreasing we have $0\leq \ell_*<\rho$. Thus, the  continuity of $n_f$ in $[0, \rho)$ and \eqref{eq:conv1} imply  $0\leq \ell_{*}=|n_f(\ell_*)|$ and from  Proposition~\ref{pr:incr3} we have $\ell_{*}=0$. Therefore,  the convergence of $\{x_k \}$ to $x_*$ is proved. The uniqueness was proved in Lemma~\ref{pr:uniq}.

For  proving   the  equality in  \eqref{eq:q2e}  note that equation \eqref{eq:conv1} implies
$$
\left[\|x_{k+1}-x_*\|\big{/}\|x_{k}-x_*\|\right]\leq \left[|n_f(\|x_{k}-x_*\|)|\big{/}\|x_{k}-x_*\|\right], \qquad k=0,1, \ldots.
$$
Since $\lim_{k\to \infty}\|x_{k}-x_*\|=0$ the  desired equality follows from first statement in  Proposition~\ref{pr:incr3}.

Now we will show  \eqref{eq:q3e}. First,  we will prove by induction that the sequences  $\{t_k \}$ and  $\{x_k \}$ defined, respectively, in \eqref{NFS} and \eqref{eq:tknk} satisfy
\begin{equation}\label{eq:mjs}
\|x_{k}-x_*\|\leq t_k, \qquad k=0,1, \ldots.
\end{equation}
Because $t_0=\|x_0-x_*\|$,  the above inequality holds for $k=0$. Now, assume that 
$\|x_{k}-x_*\|\leq t_k$.  Using \eqref{NFS},  Lemma~\ref{le:cl2}, the induction assumption  and \eqref{eq:tknk} we obtain that 
$$
\|x_{k+1}-x_*\|=\|N_F(x_k)-x_*\|\leq \frac{|n_f(t_k)|}{t_{k}^{p+1}}\,\|x_{k}-x_*\|^{p+1}\leq|n_f(t_k)|=t_{k+1},
$$
and the  proof by induction is complete.  Therefore, it easy to see that the desired inequality follows by combining  \eqref{NFS}, \eqref{eq:mjs},  Lemma~\ref{le:cl2}  and \eqref{eq:tknk}.
\end{proof}
The proof of Theorem~\ref{th:nt} follows from Corollary~\ref{cr:kanttk}, Lemmas~\ref{pr:uniq} and \ref{pr:best} and Proposition~\ref{pr:nthe}.
\section{Special Cases} \label{apl}
In this section, we will present two  special cases of Theorem~\ref{th:nt}.
\subsection{Convergence result  under H\"{o}lder-like  condition}
In this section we will present the  convergence theorem for Newton's method under affine invariant H\"{o}lder-like condition which has appeared   in   \cite{Hu04} and \cite{W03}.

\begin{theorem} \label{th:HV}
Let $\banacha$, $\banachb$ be Banach spaces, $\Omega\subseteq \banacha$ an open set and
  $F:{\Omega}\to \banachb$  a continuously
  differentiable function. 
 Let $x_*\in \Omega$ and $\kappa:=\sup\{t\in [0, R): B(x_*, t)\subset \Omega\}$.
  Suppose that $F(x_*)=0$,  $F '(x_*)$ is  invertible and there exists a constant $K>0$ and $ 0< p \leq 1$ such that  
\begin{equation} \label{eq:hc}
\left\|F'(x_*)^{-1}\left[F'(x)-F'(x_*+\tau(x-x_*))\right]\right\|\leq  K(1-\tau^p) \|x-x_*\|^p, \qquad   x\in B(x_*, \kappa) \quad \tau \in [0,1].
\end{equation}
Let $r=\min \{\kappa, \,[(p+1)/((2p+1)K)]^{1/p}\}$. Then, the sequences  with starting point $x_0\in B(x_*, r)/\{x_*\}$ and $t_0=\|x_0-x_*\|$, respectively,
$$
    x_{k+1} ={x_k}-F'(x_k) ^{-1}F(x_k), \qquad t_{k+1} =\frac{K\,p\,t_{k}^{p+1}}{(p+1)[1-K\,t_k^{p}]},\qquad k=0,1,\ldots\,, 
$$
are well defined, $\{t_k\}$ is strictly decreasing, is contained in $(0, r)$ and converges to $0$ and $\{x_k\}$ is contained in $B(x_*,r)$,  converges to  $x_*$ which is the unique zero of $F$ in $B(x_*, \,[(p+1)/K]^{1/p})$ and there holds
$$
 \|x_{k+1}-x_*\|
\leq \frac{K\,p}{(p+1)[1-K\,t_k^{p}]}\,\|x_{k}-x_*\|^{p+1}, \qquad k=0,1,\ldots\,
$$
Moreover, if  $[(p+1)/((2p+1)K)]^{1/p}<\kappa$ then $r=[(p+1)/((2p+1)K)]^{1/p}$ is the best  possible convergence radius.
\end{theorem}
\begin{proof} 
It is immediate to prove that  $F$, $x_*$ and $f:[0, \kappa)\to \mathbb{R}$, defined by 
$
f(t)=Kt^{p+1}/(p+1)-t, 
$
satisfy the inequality \eqref{Hyp:MH} and the conditions  {\bf h1}, {\bf h2} and  {\bf h3} in Theorem \ref{th:nt}. In this case, it is easy to see that  $\rho$ and $\nu$, as defined in Theorem \ref{th:nt}, satisfy
$$
\rho=[(p+1)/((2p+1)K)]^{1/p} \leq \nu=[1/K]^{1/p}, 
$$
and, as a consequence,  $r=\min \{\kappa,\; [(p+1)/((2p+1)K)]^{1/p}\}$. Moreover,   $f(\rho)/(\rho f'(\rho))-1=1$, $f(0)=f([(p+1)/K]^{1/p})=0$ and $f(t)<0$ for all $t\in (0,\, [(p+1)/K]^{1/p})$.
Therefore, the result follows  by invoking Theorem~\ref{th:nt}.
\end{proof}
\begin{remark}
Since Theorem~\ref{th:HV} is a special case of Theorem~\ref{th:nt} it follows from   Remark~\ref{r:rqc} that  
$$
\|x_{k}-x_*\| \leq  \left[ \frac{K\,p\,\|x_0-x_*\|^{p}}{(p+1)[1-K\,\|x_0-x_*\|^{p}]}\right]^{[(p+1)^k-1]/p}\,\|x_0-x_*\|, \qquad k=0,1,\ldots .
$$
\end{remark}
\begin{remark}
If  $F:{\Omega}\to \banachb$ satisfies the Lipschitz condition 
$
\|F'(x)-F'(y)\| \leq L \|x-y\|,
$
for all $ x,\,y\in \Omega,$
where $L>0$, then it also satisfies the condition \eqref{eq:hc} with $p=1$ and  $K=L\| F'(x_*)^{-1}\|$. In this case, the best  possible convergence radius for Newton's method is  $r=2/(3L\|F '(x_*)^{-1} \|)$, see~\cite{r74} and \cite{WCR77}.
We point out that the convergence radius of affine invariant  theorems are
 insensitive to invertible linear transformation of $F$. On the other hand,   theorems with the  Lipschitz condition are sensitive, see  \cite{TW79}. For more details about  affine invariant theorems on Newton's method  see \cite{DHB} (see also  \cite{DH}).
\end{remark}
\subsection{Convergence result  under generalized Lipschitz condition}
In this section, we will present a local convergence theorem on Newton's method  under a generalized Lipschitz condition due to X. Wang,  it has appeared   in  \cite{W03} (see also \cite{W00}). It is worth point out that the result in this section does not assume that  the function which defines the generalized Lipschitz condition  is  nondecreasing.

\begin{theorem} \label{th:XWT}
Let $\banacha$, $\banachb$ be Banach spaces, $\Omega\subseteq \banacha$ an open set and
  $F:{\Omega}\to \banachb$  a continuously
  differentiable function. 
 Let $x_*\in \Omega$ and  $\kappa:=\sup\{t\in [0, R): B(x_*, t)\subset \Omega\}$.
  Suppose that $F(x_*)=0$,  $F '(x_*)$ is  invertible and there exists a positive  integrable function $L:[0,\; R)\to \mathbb{R}$ such that  
\begin{equation}\label{Hyp:XW}
\left\|F'(x_*)^{-1}\left[F'(x)-F'(x_*+\tau(x-x_*))\right]\right\| \leq  \int^{\|x-x_*\|}_{\tau\|x-x_*\|} L(u){\rm d}u,
\end{equation}
for all $\tau \in [0,1]$, $x\in B(x_*, \kappa)$. Let $\bar{\nu}> 0 $  be the constant defined by
$$
\bar{\nu}:=\sup \left\{t\in [0, R): \displaystyle\int_{0}^{t}L(u){\rm d}u-1< 0\right\},
$$
and let $\bar{\rho}> 0 $ and $\bar{r}>0$  be the constants defined by
$$
\bar{\rho}:=\sup \left\{t\in (0, \delta):
\displaystyle\int^{t}_{0}L(u)u {\rm d}u\Big{/}\left[t\left(1-\displaystyle\int^{t}_{0}L(u){\rm d}u\right)\right]<1, \; t\in (0, \delta)\right\}, \qquad \bar{r}=\min \left\{\kappa, \bar{\rho}\right\}.
$$
Then, the sequences  with starting point $x_0\in B(x_*, \bar{r})/\{x_*\}$ and $t_0=\|x_0-x_*\|$, respectively,
$$
    x_{k+1} ={x_k}-F'(x_k) ^{-1}F(x_k), \qquad t_{k+1} =\displaystyle\int^{t_k}_{0}L(u)u {\rm d}u\Big{/}\left(1-\displaystyle\int^{t_k}_{0}L(u){\rm d}u\right),\qquad k=0,1,\ldots\,, 
$$
are well defined, $\{t_k\}$ is strictly decreasing, is contained in $(0, \bar{r})$ and converges to $0$, $\{x_k\}$ is contained in $B(x_*,\bar{r})$,  converges to  $x_*$ which is the unique zero of $F$ in $B(x_*, \bar{\sigma})$, where 
$$
\bar{\sigma}:=\sup\left \{t\in(0, \kappa):  \int^{t}_{0}L(u)(t-u){\rm d}u- t< 0 \right\}.
$$
and there hold: $\lim_{k\to \infty}t_{k+1}/t_k=0$ and
$ 
\lim_{k\to \infty}[\|x_{k+1}-x_*\|/\|x_k-x_*\|]=0. 
$
Moreover, if 
$$
\displaystyle\int^{\bar{\rho}}_{0}L(u)u {\rm d}u\Big{/}\left[\bar{\rho}\left(1-\displaystyle\int^{\bar{\rho}}_{0}L(u){\rm d}u\right)\right]= 1,
$$
 and ${\bar \rho}<\kappa$ then $\bar{r}=\bar{\rho}$ is the best  possible convergence radius.

\noindent
If, additionally, given $0\leq p\leq1$
\begin{itemize}
  \item[{\bf  h)}] the function  
  $  
 (0,\, \nu) \ni t \mapsto t^{1-p}L(t)
 $
is  nondecreasing,
\end{itemize}
then the sequence $\{t_{k+1}/t_k^{p+1}\}$ is strictly decreasing  and there holds
\begin{equation}
    \label{eq:wq3}
\|x_{k+1}-x_*\| \leq \big[t_{k+1}/t_k^{p+1}\big]\,\|x_k-x_*\|^{p+1}, \qquad k=0,1,\ldots\,. 
  \end{equation}
\end{theorem}
\begin{proof} 
Let  ${\bar f}:[0, \kappa)\to \mathbb{R}$ a differentiable function defined by 
\begin{equation} \label{eq:wf}
{\bar f}(t)=\int_{0}^{t}L(u)(t-u){\rm d}u-t. 
\end{equation}
Note  that the derivative of the function $f$ is given by 
$$
 {\bar f}'(t)=\int_{0}^{t}L(u){\rm d}u-1.
$$
Because $L$ is integrable ${\bar f}'$ is continuous (in fact ${\bar f}'$ is absolutely continuous). So, it is easy to see that \eqref{Hyp:XW} becomes \eqref{Hyp:MH} with $f'={\bar f}'$. Moreover, because $L$ is positive the function $f={\bar f}$  satisfies  the conditions  {\bf h1} and  {\bf h2} in Theorem \ref{th:nt}. Direct algebraic manipulation yields
$$
\frac{1}{t^{p+1}} \left[\frac{{\bar f}(t)}{{\bar f}'(t)}-t\right]=\left[
\frac{1}{t^{p+1}}\displaystyle\int^{t}_{0}L(u)u{\rm d}u\right]
\frac{1}{|{\bar f}'(t)|}.
$$
If assumption {\bf  h} holds then Lemma~$2.2$ of \cite{W03} implies that the first tern in right had side of the above equation is nondecreasing  in $(0,\, \nu)$. Now,  since $1/|{\bar f}'|$ is  strictly  increasing in $(0,\, \nu)$ the above equality implies that {\bf  h3} in Theorem~\ref{th:nt},  with $f={\bar f}$, also holds.
Therefore, the result  follows from Theorem~\ref{th:nt} with $f={\bar f}$, $\nu=\bar{\nu}$, $\rho=\bar{\rho}$, $r=\bar{r}$ and $\sigma=\bar{\sigma}$.
\end{proof}
\begin{remark}
Since Theorem~\ref{th:XWT} is a special case of Theorem~\ref{th:nt} it follows from   Remark~\ref{r:rqc} that if  $p=0$ then $\|x_{k}-x_*\| \leq q^{k}\,\|x_0-x_*\|$, for $k=0,1,\ldots$ and if  $0<p\leq 1$ then
$$
\|x_{k}-x_*\| \leq q^{[(p+1)^k-1]/p}\,\|x_0-x_*\|, \qquad k=0,1,\ldots ,
$$
where
$$
q=\displaystyle\int^{\|x_0-x_*\|}_{0}L(u)u {\rm d}u\Big{/}\left[\|x_0-x_*\|\left(1-\displaystyle\int^{\|x_0-x_*\|}_{0}L(u){\rm d}u\right)\right].
$$
\end{remark}
\begin{remark}
It was shown in  \cite{W00} that if $L$ is positive and  nondecreasing then  the sequence generated by Newton's method converges with quadratic rate.  From Theorem~\ref{th:XWT} we conclude that the assumption on the nondecrement   of $L$ is needed only to obtain the quadratic convergence rate of the  Newton's sequence.  This result was also obtained in \cite{W03}. 

Finally, we observe that if the positive  integrable function $L:[0,\; R)\to \mathbb{R}$ is   nondecreasing then the strictly increasing function $f':[0,\; R)\to \mathbb{R}$, defined by
$$
 f'(t)=\int_{0}^{t}L(u){\rm d}u-1,
$$
is   convex. In this case, is not hard to prove that the inequalities \eqref{Hyp:MH} and \eqref{Hyp:XW} are equivalents. On the other hand,  if $f'$ is  strictly increasing and non necessary convex  then the inequalities  \eqref{Hyp:MH} and \eqref{Hyp:XW} are not equivalents. Because there exists functions strictly increasing, continuous, with  derivative zero almost everywhere, see \cite{T78} (see also \cite{OW07}). Note  that these functions  are not absolutely continuous, so they can not be  represented by an integral.
\end{remark}
\section{Final remarks } \label{fr}
Theorem \ref{th:XWT} has  many interesting special cases, including the  Smale's theorem  on Newton's method (see \cite{S86}) for analytical functions, see \cite{W03}.  Theorem \ref{th:nt} has   the  Nesterov-Nemirovskii's theorem on Newton's method (see \cite{NN-94}) for self-concordant functions  as a special case, see \cite{F08}. 

\end{document}